\documentclass[11pt,reqno]{article}
\usepackage{jmlr2e}
\usepackage{fullpage,times,graphicx,amssymb,amsmath,amsfonts,bbm,psfrag,xcolor}

\usepackage{graphicx}
\usepackage{amssymb}
\usepackage{amsmath}
\usepackage{amsfonts}
\usepackage{bbm}
\usepackage{psfrag}
\usepackage{xcolor}
\usepackage{fullpage}
\usepackage{epstopdf}
\usepackage{graphicx}
\usepackage{times}
\usepackage{jmlr2e}
\usepackage{commath}
\usepackage{thmtools}
\usepackage{thm-restate}
\usepackage{times}
\usepackage{graphicx}
\usepackage{natbib}
\usepackage{dirtytalk}
\definecolor{ddarkbrown}{rgb}{0.5,0.2,0.05} \definecolor{bbluegray}{rgb}{0.05,0,0.5}

\usepackage[colorlinks,citecolor=bbluegray,linkcolor=ddarkbrown,urlcolor=blue,breaklinks]{hyperref}

\usepackage{subfig,float} 

\usepackage{mathtools}

\usepackage{algorithm,algcompatible,algpseudocode}
\usepackage{multicol,lipsum}
\algnewcommand{\Inputs}[1]{%
	\State \textbf{Inputs: \:}{#1}
}

\algnewcommand{\Output}[1]{%
	\State \textbf{Output: \:}{#1}
}
\algnewcommand{\Initialize}[1]{%
	\State \textbf{Initialize: \:}{#1}
}

\algnewcommand{\IIf}[1]{\State\algorithmicif\ #1\ \algorithmicthen}
\algnewcommand{\EndIIf}{\unskip\ \algorithmicend\ \algorithmicif}

\let \oldsection \section
\renewcommand{\section}{\vspace{3ex plus 1ex}\oldsection}

\newcommand{\BEAS}{\begin{eqnarray*}}
	\newcommand{\EEAS}{\end{eqnarray*}}
\newcommand{\BEA}{\begin{eqnarray}}
\newcommand{\EEA}{\end{eqnarray}}
\newcommand{\mb}{\mathbb}
\newcommand{\BEQ}{\begin{equation}}
\newcommand{\EEQ}{\end{equation}}
\newcommand{\BIT}{\begin{itemize}}
	\newcommand{\EIT}{\end{itemize}}
\newcommand{\BNUM}{\begin{enumerate}}
	\newcommand{\ENUM}{\end{enumerate}}
	\newcommand{\la}{\langle}
	\newcommand{\ra}{\rangle}

	\newcommand{\mr}{\mathrm}

	\newcommand{\E}{\mathbb{E}}	
	\newcommand{\hy}{\hat{y}}
	\newcommand{\by}{\Bar{y}}

	\newcommand{\iT}{{T_k^{-1}}}
	\newcommand{\iS}{{S_k^{-1}}}

\newcommand{\BA}{\begin{array}}
	\newcommand{\EA}{\end{array}}

 \numberwithin{dummy}{section}

\numberwithin{mythm}{section}
\numberwithin{mydef}{section}
\numberwithin{myprop}{section}
\numberwithin{mylem}{section}
\numberwithin{mycor}{section}

\title{Stochastic Primal-Dual Three Operator Splitting Algorithm with Extension to Equivariant  Regularization-by-Denoising}

\begin{document}
	\author{\name Junqi Tang, Matthias Ehrhardt,  Carola-Bibiane Sch\"onlieb  \email j.tang.2@bham.ac.uk\\
		\addr School of Mathematics, University of Birmingham\\ Department of Mathematical Sciences, University of Bath\\DAMTP, University of Cambridge
		}
	\editor{}
	
	
	\maketitle


\begin{abstract}
 In this work we propose a stochastic primal-dual three-operator splitting algorithm (TOS-SPDHG) for solving a class of convex three-composite optimization problems. Our proposed scheme is a direct three-operator splitting extension of the SPDHG algorithm [Chambolle et al. 2018]. We provide theoretical convergence analysis showing ergodic $O(1/K)$ convergence rate, and demonstrate the effectiveness of our approach in imaging inverse problems. Moreover, we further propose TOS-SPDHG-RED and TOS-SPDHG-eRED which utilizes the regularization-by-denoising (RED) framework to leverage pretrained deep denoising networks as image priors for improved reconstruction.

\end{abstract}
\section{Introduction}

In many imaging inverse problems and machine learning applications, we often need to solve convex composite optimization problems of the form:
\begin{align}\label{obj}
   x^\star = \min_{x \in X} f(Ax) + g(x) + h(x)
\end{align}
where $X \subseteq \mathbb{R}^d$ is a convex set, $A \in \mb{R}^{N \times d}$ is a linear operator, $f(Ax)$ the data-fidelity, $g(x)$ a regularizer which admits simple proximal operator, and $h(x)$ another regularizer for which we have access to its gradients. Here we assume that both three terms are proper convex lower-semicontinuous functions, while $h(x)$ has $L$-Lipschitz continuous gradients. Since this type of optimization problems are usually large-scale and high-dimensional for modern machine learning and imaging tasks, fast and accurate iterative algorithms are highly desirable in real world applications \citep{xiao2014proximal,defazio2014saga,allen2017katyusha, chambolle2024stochastic}. The generic composite optimization problem \eqref{obj} requires three operator splitting schemes as the optimizer \citep{raguet2013generalized,davis2017three}. Since for most applications the linear operator $A$ is non-trivial, for the case where $f$ is non-smooth (or have a large gradient Lipschitz constant), it is often more beneficial to dualize $f(Ax)$ and reformulate the objective as the saddle-point problem:
\begin{equation}\label{saddle}
  [x^\star, y^\star] =  \min_{x \in X} \max_{y \in Y} \Phi(x,y):= h(x) + g(x) + \la A x, y\ra - f^*(y),
\end{equation}
where:
\begin{equation}
    f^*(y) = \sup_{z \in Y} \la z, y\ra - f(z).
\end{equation}
A number of primal-dual first-order algorithms have been proposed for solving the saddle point problem \eqref{saddle}, and among this line of works, the most noticeable one would be the Condat-Vu's scheme \citep{condat2013primal,vu2013splitting}. In this work, based on the framework of the stochastic primal-dual hybrid gradient (SPDHG) by \cite{chambolle2018stochastic}, we propose a stochastic gradient extension of the Condat-Vu algorithm for efficiently solving the saddle point optimization problem \eqref{saddle} if it admits finite-sum structure which is true for most applications in machine learning and computational imaging.

\section{Stochastic primal-dual hybrid gradient and three operator splitting}
We start by recalling the stochastic primal-dual hybrid gradient (SPDHG) algorithm proposed by \cite{chambolle2018stochastic} for two-composite problems. For simplicity we illustrate here their uniform serial sampling $p_i = p = 1/n$ case:
\begin{align}
    f(y) = \sum_{i=1}^n f_i(y_i), \ \ 
    (Ax)_i = A_i x,\ \ A_i \in \mb{R}^{m \times d}.
\end{align}
Now we can formulate the objective to the finite-sum saddle-point form:
\begin{equation}
  [x^\star, y^\star] =  \min_{x \in X} \max_{y \in Y} \Phi(x,y):= h(x) + g(x) + \sum_{i=1}^n\la A_i x, y_i \ra - f_i^*(y_i).
\end{equation}
The most basic version of SPDHG algorithm \citep{chambolle2018stochastic} for the case $h(x)=0$ can be written as (suppose we randomly select an index $j_k$ at iteration $k$, and we denote here primal and dual step-sizes to be $\tau$ and $\sigma$ respectively):
\begin{align}
    x^{k+1} &= \operatorname{prox}_{g}^\tau(x^k - \tau z^k) \\
    y^{k+1}_i &= \begin{cases}
    \operatorname{prox}_{f^*_i}^\sigma (y^k_i + \sigma A_i x^{k+1}) & \text{if $i = j_k$} \\
    y^k_i
    \end{cases} \\
    z^{k+1} &= A^* y^{k+1} + p^{-1} A^*_{j_k} (y^{k+1}_i - y^k_i) = A^* y^{k} + (p^{-1}+1) A^*_{j_k} (y^{k+1}_i - y^k_i)
\end{align}
where the proximal operator is defined as:
\begin{equation}
    \operatorname{prox}_g^\tau(x) := \arg\min_{z \in C} \frac{1}{2}\|z - x\|_{\tau^{-1}}^2 + g(z)
\end{equation}
where the weighted norm $\|x\|_{\tau^{-1}}^2 = \la\tau^{-1}x, x\ra$ and $C$ being some constraint set. In this work, following the notation in \cite{chambolle2018stochastic,ehrhardt2017faster} we mostly use a generalized matrix form of proximal operator that allows variable metric step sizes: \[\operatorname{prox}_g^T(x) := \arg\min_{z \in C} \frac{1}{2}\|z - x\|_{T}^2 + g(z),\] where the step size matrix $T$ is symmetric positive definite, and $\|x\|_T^2 = \langle T^{-1}x, x\rangle$.

Now for the case where we additionally have a smooth regularization term $h(x)$ for which we can compute the gradient but not necessarily the proximal operator, we propose TOS-SPDHG (Three operator splitting with stochastic primal-dual hybrid gradient), which is a simple extension of the SPDHG algorithm:


 \begin{eqnarray*}
 && \mathrm{\textbf{TOS-SPDHG}} - \mathrm{Initialize}\ x^0\in \mb{R}^d, \ y^0 \in \mb{R}^q,\ \by = y^0\\
 &&\mathrm{For} \ \ \ k = 0, 1, 2,...,  K\\
&&\left\lfloor
\begin{array}{l}
    \mr{Select} \ j_k \in \{1,...,n\}\\
    x^{k+1} = \mathrm{prox}_{g}^{T_k}(x^k - T_k (A^T\by^k + \nabla h(x^k))) \\
    y^{k+1}_i = \begin{cases}
    \operatorname{prox}_{f^*_i}^{S_k^i}(y^k_i + S_k^i A_i x^{k+1}) \ \ \ \ \text{if $i = j_k$} \\
    y^k_i
    \end{cases} \\
    \by^{k+1} = y^{k+1} + \theta Q(y^{k+1} - y^k)
\end{array}
\right.
 \end{eqnarray*}

Here we denote $S_k=\mr{diag}(S_k^1,...,S_k^n)$ and $T_k$ the step size matrices\footnote{our current analysis in this note restricts only to the case of fixed step size matrices $S_k = S$, $T_k = T$ but in our ongoing works we found in practice that these can be varying and adaptively chosen for faster rates. }, and scalar $\theta > 0$ while $Q:= \mr{diag}(p_1^{-1}I,..., p_n^{-1}I)$ where each of the $p_i$ denotes the probability $i$-th index is sampled in each iteration. Although this extension of SPDHG is simple and natural, we believe that it is a crucial extension especially for imaging inverse problems. For example, for Positron Emission Tomography imaging tasks, practitioners typically use jointly the  Kullback-Leibler (KL) data-fit and the relative difference prior (the Q.Clear framework by GE Healthcare) which is a smooth regularizer with simple gradient evaluation but inefficient proximal evaluation, to cope with the extremely noisy measurement data \cite{ehrhardt2017faster}. In such case the original SPDHG is not directly applicable and we need this TOS-SPDHG extension.

Noting that just as the original SPDHG, we keep the non-scalar step-sizes $T$, $S$ for TOS-SPDHG which can be chosen as preconditioners for accelerating empirical convergence \citep{ehrhardt2017faster}. Our convergence theory is built upon the SPDHG analysis by \citep{chambolle2018stochastic} and hence we also support arbitrary sampling schemes, for example we might use importance sampling to accelerate convergence in some applications.

\clearpage

\subsection{Extension to Equivariant RED Priors}

Our proposed TOS-SPDHG can easily incorporate with advanced image prior such as the regularization-by-denoising (RED) \citep{romano2017little} with denoising algorithms (BM3D/NLM) and deep denoising networks (DnCNN). In this work, we propose two extensions for TOS-SPDHG with RED and recently introduced equivariant denoisers \citep{terris2024equivariant}.

Our first extension is the TOS-SPDHG-RED algorithm which we directly include RED-step inside our TOS-SPDHG, which allow us to utilize jointly the smooth prior $h(x)$ and the deep denoising prior implicit given by $x - \mathcal{D}(x)$ in the algorithm. Here we denote $\mathcal{D}$ as the denoising function, which in our experiments we choose the standard DnCNN \cite{zhang2017beyond} denoiser.

 \begin{eqnarray*}
 && \mathrm{\textbf{TOS-SPDHG-RED}} - \mathrm{Initialize}\ x^0\in \mb{R}^d, \ y^0 \in \mb{R}^q,\ \by = y^0\\
 &&\mathrm{For} \ \ \ k = 0, 1, 2,...,  K\\
&&\left\lfloor
\begin{array}{l}
    \mr{Select} \ j_k \in \{1,...,n\}\\
    x^{k+1} = \mathrm{prox}_{g}^{T_k}(x^k - T_k (A^T\by^k + \nabla h(x^k) + \lambda(x^k - \mathcal{D}(x^k)))) \\
    y^{k+1}_i = \begin{cases}
    \operatorname{prox}_{f^*_i}^{S_k^i}(y^k_i + S_k^i A_i x^{k+1}) \ \ \ \ \text{if $i = j_k$} \\
    y^k_i
    \end{cases} \\
    \by^{k+1} = y^{k+1} + \theta Q(y^{k+1} - y^k)
\end{array}
\right.
 \end{eqnarray*}

Moreover, we further consider recently introduced equivariant-denoiser trick \cite{terris2024equivariant} for improved numerical results and stability while using advanced deep denoising networks such as DnCNN and DRUNet. Denote unitary matrices $\{\mathcal{T}_{r_i}\}_{i=1}^{|G|}$ as the transforms for some group $G$ (for instance, rotation for tomographic image reconstruction), we can write the corresponding denoising step for RED as

\begin{equation}
    x^k - \mathcal{T}_{r_k}^{-1}\mathcal{D}(\mathcal{T}_{r_k} x^k),\ \ \ r_k \sim G
\end{equation}
 Plugging this into TOS-SPDHG we obtain our new algorithm TOS-SPDHG-eRED:
  \begin{eqnarray*}
 && \mathrm{\textbf{TOS-SPDHG-eRED}} - \mathrm{Initialize}\ x^0\in \mb{R}^d, \ y^0 \in \mb{R}^q,\ \by = y^0\\
 &&\mathrm{For} \ \ \ k = 0, 1, 2,...,  K\\
&&\left\lfloor
\begin{array}{l}
    \mr{Select} \ j_k \in \{1,...,n\}, \ r_k\sim G\\
    x^{k+1} = \mathrm{prox}_{g}^{T_k}(x^k - T_k (A^T\by^k + \nabla h(x^k) + \lambda(x^k - \mathcal{T}_{r_k}^{-1}\mathcal{D}(\mathcal{T}_{r_k} x^k)))) \\
    y^{k+1}_i = \begin{cases}
    \operatorname{prox}_{f^*_i}^{S_k^i}(y^k_i + S_k^i A_i x^{k+1}) \ \ \ \ \text{if $i = j_k$} \\
    y^k_i
    \end{cases} \\
    \by^{k+1} = y^{k+1} + \theta Q(y^{k+1} - y^k)
\end{array}
\right.
 \end{eqnarray*}
We numerically observe that the proposed RED extensions can be effectively applied using advanced deep denoising networks, especially the later TOS-SPDHG-eRED scheme.

\section{Convergence Analysis for TOS-SPDHG}

We start by introducing the definition of expected separable overapproximation (ESO) property \citep{richtarik2014iteration} which is standard for the analysis of coordinate descent \citep{fercoq2015accelerated}, stochastic variance-reduced gradient methods \citep{konevcny2016mini}, and stochastic primal-dual gradient methods \citep{chambolle2018stochastic} with arbitrary sampling:
\begin{definition}
\citep{richtarik2014iteration} (Expected Separable Overapproximation) Let $\mb{S} \in \{1,...,n\}$ and $p_i := \mb{P}(i\in \mb{S})$, and bounded linear operators $C_i$ and $(Cx)_i:= C_i x$, then the ESO parameters $\{v_i\}$ satisfy:
\begin{equation}
    \mb{E}_{\mb{S}}\|\sum_{i \in \mb{S}} C_i^Tz_i\|_2^2 \leq \sum_{i=1}^n p_i v_i \|z_i\|_2^2
\end{equation}
\end{definition}
As we will see, the choice of $T,S$ for provable convergence of TOS-SPDHG depends on the ESO parameters of the linear operator $S^{\frac{1}{2}}A(T^{-1} - LI)^{-\frac{1}{2}}$, where $L$ is the gradient Lipschitz constant of $h$.

\subsection{Convergence analysis}

For the convergence proof of TOS-SPDHG algorithm we need following two lemmas.
\begin{lemma}
Assuming that $h$ has $L$-Lipschitz continuous gradients, and both $f$, $g$, and $h$ are proper convex lower-semicontinous functions, for the deterministic updates
\begin{equation}
\begin{aligned}
&x^{k+1} = \mathrm{prox}_g^{T_k}(x^k - T_k (A^T\Bar{y}^k + \nabla h(x^k)))\\
&\hat{y}_i^{k+1} = \mathrm{prox}_{f_i^*}^{S_k^i} (y_i^k + S_k^i A_i x^{k+1}), \ \ i = 1,...,n,
\end{aligned}
\end{equation}
we have:
\begin{equation}
\begin{aligned}
&\|x^k - x\|_{T_k^{-1}}^2 + \|y^k - y\|_{S_k^{-1}}^2\\
\geq& \|x^{k+1} - x\|_{T_k^{-1} } + \|\hy^{k+1} - y\|_{S_k^{-1} } \\
& + 2\left(g(x^{k+1})-g(x)+h(x^{k+1})-h(x)+f^*(\hy^{k+1})-f^*(y)\right) - 2\langle A(x - x^{k+1}), \by^k \rangle\\
&+ \|x^{k+1} - x^k\|_{T_k^{-1} -LI}^2 + \|y^{k+1} - y^k\|_{S_k^{-1}}^2\\
\end{aligned}
\end{equation}

\end{lemma}
We include the proof in the first half of next subsection. The following Lemma can be directly derived from \citep[Lemma 4.2]{chambolle2018stochastic}:
\begin{lemma}
Let $\mb{S} \subseteq \{1,...,n\}$ be a random set, $y_i^+ = \hat{y}_i$ if $i \in \mb{S}$ and $y_i$ otherwise, and $\{v_i\}$ be some ESO parameter of $S^{\frac{1}{2}}A(T^{-1} - LI)^{-\frac{1}{2}}$, then for any $c>0$ we have:
\begin{equation}
    2\mb{E}_{\mb{S}}\langle QAx, y^+ - y \rangle \geq -\mb{E_{\mb{S}}}\left\{ \frac{1}{c}\|x\|_{T^{-1} - LI} + c \max_i \frac{v_i}{p_i}\|y^+ - y\|_{QS^{-1}}^2\right\}
\end{equation}
\end{lemma}

Now following similar steps of the proof for \cite[Theorem 4.3]{chambolle2018stochastic}, with a different step-size choice taking into account the $L$-smoothness of $h$, we can have the following convergence guarantee for TOS-SPDHG:
\begin{theorem}
Assuming that $h$ has $L$-Lipschitz continuous gradients, and both $f$, $g$, and $h$ are proper convex lower-semicontinous functions, let $\theta=1$ and fixed step-size matrices $T, S$ be selected such that $(T^{-1} - LI)$ is positive-definite and there exist ESO parameters $\{v_i\}$ of $S^{\frac{1}{2}}A(T^{-1} - LI)^{-\frac{1}{2}}$ with $v_i < p_i, \forall i \in [n]$, then we have:
\begin{equation}
    \mathbb{E}\left(\Phi(x_{(K)}, y^\star) - \Phi(x^\star, y_{(K)})\right) \leq O\left( \frac{1}{K} \right)
\end{equation}
where $x_{(K)} = \frac{1}{K}\sum_{k=1}^K x^K$, $y_{(K)} = \frac{1}{K}\sum_{k=1}^K y^K$ denotes the ergodic sequences.
\end{theorem}

\subsection{Convergence proof}

The proof of the main theorem follows similar steps in \citep{chambolle2018stochastic} which was originally for two-composite problems, while here we extend their framework to three-composite setting. The first part of the proof is regarding Lemma 1. By the definition of proximal operator we have:
\begin{equation}
    x^{k+1} - x^k + T_k(A^T\by^{k} + \nabla h(x^k)) \in T_k \partial g(x^{k+1})
\end{equation}
and due to the convexity of $g$, for any $x$ we have:
\begin{equation}
    g(x) - g(x^{k+1}) - \langle \partial g(x^{k+1}), x - x^{k+1} \rangle \geq 0,
\end{equation}
and hence:
\begin{equation}
\begin{aligned}
&g(x) \geq g(x^{k+1}) + \langle T_k^{-1} (x^k - x^{k+1}) - A^T \by^k - \nabla h(x^k), x - x^{k+1} \rangle  \\
\end{aligned}
\end{equation}
Identically, by the definition of proximal operator and convexity of $f_i^*$, for any $y$ we can have:
\begin{equation}
\begin{aligned}
&f_i^*(y_i) \geq f_i^*(\hy_i^{k+1}) + \langle (S_k^i)^{-1} (y_i^k - \hy_i^{k+1}) + A_i x^{k+1}, y_i - \hy^{k+1} \rangle\\
\end{aligned}
\end{equation}
Summing up twice the inequalities we have:
\begin{equation}
\begin{aligned}
2f^*(y) +2g(x) &\geq 2f^*(\hy^{k+1}) + 2g(x^{k+1}) + \|x^k - x^{k+1}\|_{T_k^{-1}}^2 + \|x^{k+1} - x\|_{\iT}^2 \\&- \|x^k - x\|_\iT^2 -  \la A(x - x^{k+1}), \by^k \ra - \la \nabla h(x^k), x - x^{k+1} \ra\\
& + \|y^k - \hy^{k+1}\|_\iS^2 + \|y - \hy^{k+1}\|_{\iS}^2 - \|y^k - y\|_\iS^2\\
& + \la Ax^{k+1}, y - \hy^{k+1} \ra.
\end{aligned}
\end{equation}
(Utilizing the identity that $2\langle S(a-b), c-b\rangle = \|a-b\|_S^2 + \|b-c\|_S^2 -\|a-c\|_S^2$) Then we arrange the terms to be:
\begin{equation}
\begin{aligned}
\|x^k - x\|_\iT^2 + \|y^k - y\|_\iS^2 &\geq \|x^{k+1} - x\|_{\iT}^2 +  \|y - \hy^{k+1}\|_{\iS}^2\\
& + 2(g(x^{k+1}) - g(x) + f^*(\hy^{k+1}) - f^\star(y))\\
& + 2(\la Ax^{k+1}, y-\hy^{k+1} \ra - \la A(x-x^{k+1}), \by^k \ra) \\
& + \|y^k - \hy^{k+1}\|_\iS^2 + \|x^k - x^{k+1}\|_{T_k^{-1}}^2\\
& - 2\la \nabla h(x^k), x - x^{k+1} \ra.
\end{aligned}
\end{equation}
Noting that we have made the assumption that $h$ is $L$-smooth, hence for this extra term $- \la \nabla h(x^k), x - x^{k+1} \ra$ we can bound it using the three-point descent identity for convex and smooth functions:
\begin{equation}
     \la \nabla h(x^k), x^{k+1} - x \ra \geq h(x^{k+1}) - h(x) - \frac{L}{2}\|x^{k+1} - x^k\|_2^2,
\end{equation}
then by summing these two inequalities we finish the proof of Lemma 2:
\begin{equation}
\begin{aligned}
&\|x^k - x\|_{T_k^{-1}}^2 + \|y^k - y\|_{S_k^{-1}}^2\\
\geq& \|x^{k+1} - x\|_{T_k^{-1} } + \|\hy^{k+1} - y\|_{S_k^{-1} } \\
& + 2\left(g(x^{k+1})-g(x)+h(x^{k+1})-h(x)+f^*(\hy^{k+1})-f^*(y)\right) - 2\langle A(x - x^{k+1}), \by^k \rangle\\
&+2\la Ax^{k+1}, y-\hy^{k+1} \ra + \|x^{k+1} - x^k\|_{T_k^{-1} -LI}^2 + \|\hy^{k+1} - y^k\|_{S_k^{-1}}^2\\
\end{aligned}
\end{equation}
Now following similar steps from the proof of \cite[Theorem 4.3]{chambolle2018stochastic} let's denote $w=(x,y)$, $w'=(x',y')$ and keep the step-size matrices $T,S$ fixed $\forall k$, denote:
\begin{equation}
    F_i(y_i | x',y_i') := f_i^*(y_i) - f_i^*(y_i') - \la A_ix', y_i -y_i' \ra , \ \ \  F(y|w') := \sum_{i=1}^n F_i(y_i|x', y_i'),
\end{equation}
and:
\begin{equation}
    F^p(y|w'):= \sum_{i=1}^n (1/p_i - 1) F_i(y_i|x',y_i'),
\end{equation}
and:
\begin{equation}
    G(x|w'):= g(x)+h(x) - g(x')-h(x') - \la -A^Ty', x-x' \ra,
\end{equation}
and:
\begin{equation}
    H(w|w'):= G(x|w') + F(y|w')
\end{equation}
Then we can have:
\begin{equation}
\begin{aligned}
&\|x^k - x\|_{T^{-1}}^2 + \|y^k - y\|_{QS^{-1}}^2 + 2F^p(y^k|w)\\
\geq& \E\{\|x^{k+1} - x\|_{T^{-1}}^2 + \|y^{k+1} - y\|_{S^{-1}}^2 + 2F^p(y^{k+1}|w) +2H(w^{k+1}|w) \\
&  - 2\langle A(x^{k+1} - x), Q(y^{k+1} - y^k) + y^k - \by^k \rangle \\&+ \|x^{k+1} - x^k\|_{T^{-1} -LI}^2 + \|y^{k+1} - y^k\|_{QS^{-1}}^2\} \\
\end{aligned}
\end{equation}
Meanwhile denote $\beta:=\max_{i \in [n]} \frac{v_i}{p_i}$, by ESO we can have:
\begin{equation}
    2\mb{E}\langle QA(x^{k+1} - x^k), y^k - y^{k-1} \rangle \geq -\mb{E}\left\{ \beta\|x^{k+1} - x^k\|_{T^{-1} - LI} + \|y^k - y^{k-1}\|_{QS^{-1}}^2\right\}.
\end{equation}
Then:
\begin{equation}
\begin{aligned}
&\|x^k - x\|_{T^{-1}}^2 + \|y^k - y\|_{QS^{-1}}^2 + \|y^k - y^{k-1}\|_{QS^{-1}}^2 + 2F^p(y^k|w) \\&- 2\la QA(x^{k} - x), y^{k} -y^{k-1} \ra\\
\geq& \E\{\|x^{k+1} - x\|_{T^{-1}}^2 + \|y^{k+1} - y\|_{S^{-1}}^2 + 2F^p(y^{k+1}|w) +2H(w^{k+1}|w) \\
& - 2\la QA(x^{k+1} - x), y^{k+1} -y^{k} \ra  + 2\langle QA(x^{k+1} - x^k), y^{k} - y^{k-1} \rangle \\&+ \|x^{k+1} - x^k\|_{T^{-1} -LI}^2 + \|y^{k+1} - y^k\|_{QS^{-1}}^2\} \\
\geq& \E\{\|x^{k+1} - x\|_{T^{-1}}^2 + \|y^{k+1} - y\|_{QS^{-1}}^2 + 2F^p(y^{k+1}|w) +2H(w^{k+1}|w) \\
& - 2\la QA(x^{k+1} - x), y^{k+1} -y^{k} \ra  - \beta\|x^{k+1} - x^k\|_{T^{-1} - LI} \\
& + \|x^{k+1} - x^k\|_{T^{-1} -LI}^2 + \|y^{k+1} - y^k\|_{QS^{-1}}^2\} \\
\end{aligned}
\end{equation}
Let's denote:
\begin{equation}
\begin{aligned}
       & \Omega^k := \E\{ \|x^k - x\|_{T^{-1}}^2 + \|y^k - y\|_{QS^{-1}}^2 \\&+ \|y^k - y^{k-1}\|_{QS^{-1}}^2 + 2F^p(y^k|w) - 2\la QA(x^{k} - x), y^{k} -y^{k-1} \ra\},
\end{aligned}
\end{equation}
 for $\sqrt{\beta} < 1$ and $T^{-1} - LI$ being positive definite we have:
\begin{equation}
\begin{aligned}
\Omega^k &\geq \Omega^{k+1} + \E(2H(w^{k+1}|w) + (1-\beta)\|x^{k+1} - x^k\|_{T^{-1} - LI}) \\&\geq  \Omega^{k+1} + 2\E H(w^{k+1}|w).
\end{aligned}
\end{equation}
By summing over above inequality from $k=0$ to $K-1$, at the same time notice that $\Omega^k \geq 2 \E F^p(y^k | w)$ due to ESO, we have:
\begin{equation}
    \E \sum_{k=1}^K H(w^k|w) \leq \frac{\Omega^0}{2} - \E F^p(y^K|w) < + \infty,
\end{equation}
then due to convexity of $H$, for the ergodic sequence $w_(k) = \frac{1}{K}\sum_{k=1}^K w_k$, we have:
\begin{equation}
   \mathbb{E}\left(\Phi(x_{(K)}, y^\star) - \Phi(x^\star, y_{(K)})\right) \leq \E H(w_{(K)} | w) \leq \frac{1}{K} \E \sum_{k=1}^K H(w^k|w) = O(1/K)
\end{equation}
Thus finishes the proof of the Theorem 1.

\begin{figure}[t]
   \centering
    {\includegraphics[width= .995\textwidth]{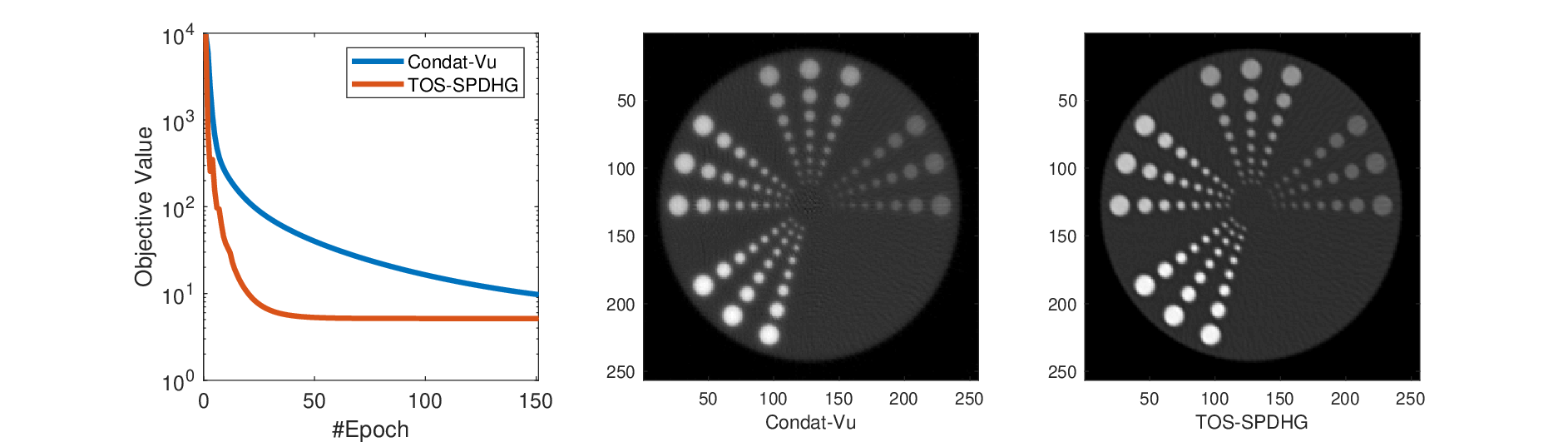}}
    {\includegraphics[width= .995\textwidth]{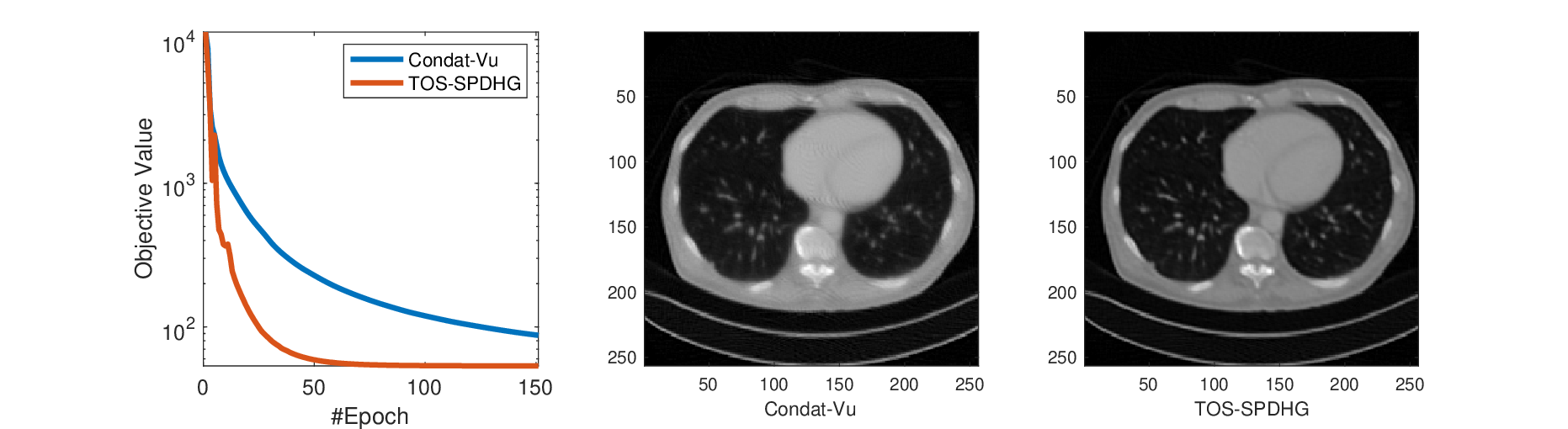}}

   \caption{Sparse-view X-ray CT Reconstruction results for Condat-Vu and TOS-SPDHG, termintate at 150th epoch. For the first example we use least-squares as data-fit, where for the second example we use the KL divergence.}
\label{F1}
\end{figure}

\begin{figure}[t]
   \centering
    {\includegraphics[width= .875\textwidth]{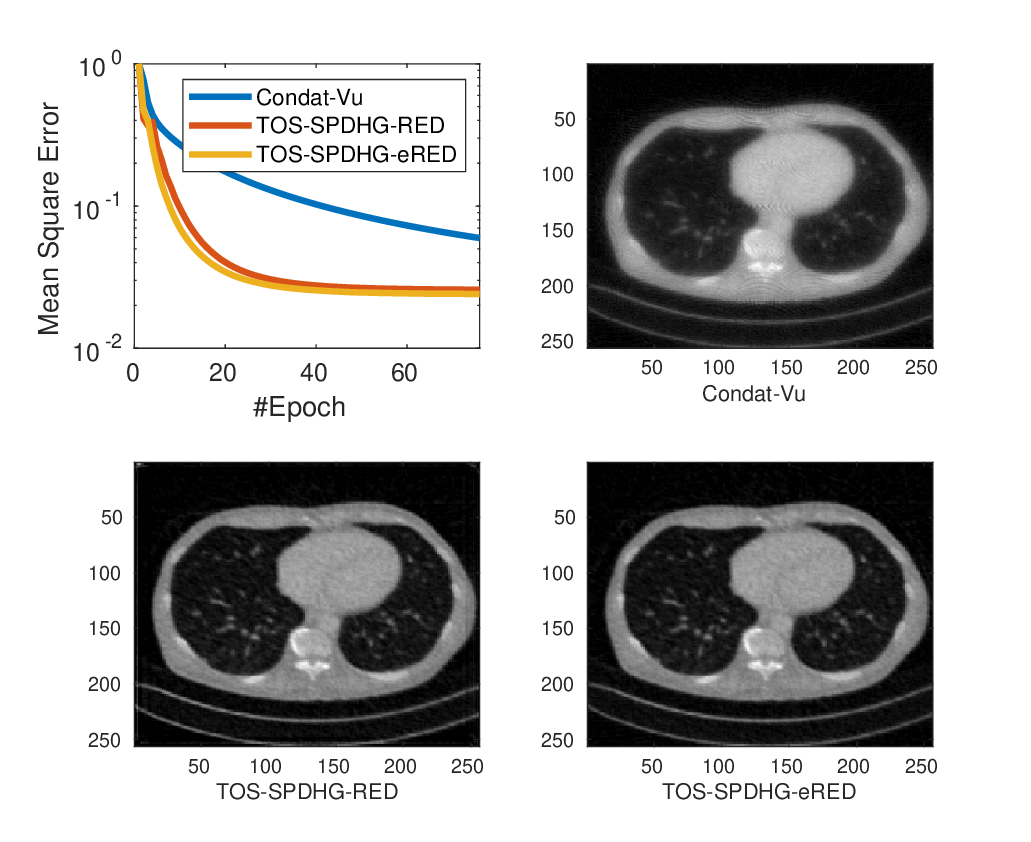}}

   \caption{Low-dose sparse-view X-ray CT Reconstruction results for Condat-Vu, TOS-SPDHG-RED and TOS-SPDHG-eRED, termintate at 75th epoch. We choose DnCNN as the denoiser.}
\label{F2}
\end{figure}

\begin{figure}[t]
   \centering
    {\includegraphics[width= .875\textwidth]{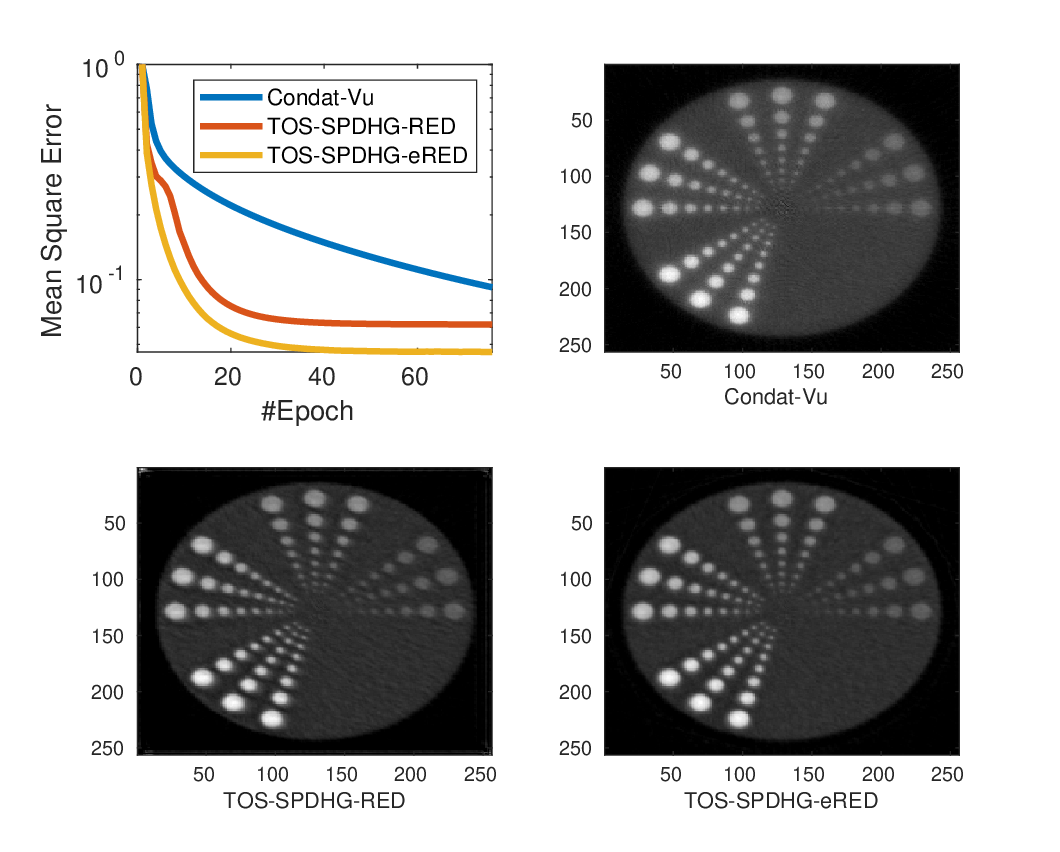}}

   \caption{Low-dose sparse-view X-ray CT Reconstruction results for Condat-Vu, TOS-SPDHG-RED and TOS-SPDHG-eRED, termintate at 75th epoch. We choose DnCNN as the denoiser.}
\label{F3}
\end{figure}

\section{Numerical Experiment}

In this section we present our preliminary numerical experiments on X-ray Computed Tomography (CT) image reconstruction which is an imaging inverse problem suitable for stochastic first-order methods \citep{tang2020practicality}, where we seek to find an estimate of some ground-truth image $x^\dagger$ via measurement $b$ corrupted by Poisson noise:
\begin{equation}
    b \sim \mathrm{Poisson}(I_0 e^{-Ax^\dagger}),
\end{equation}
via solving:
\begin{align}\label{obj1}
   x^\star = \min_{x \in X} f(Ax, b) + g(x) + h(x).
\end{align}
Here we choose the smooth edge-preserving regularizer \citep{thibault2007three} $h(x) = \lambda \sum_{d_i} \phi([Dx]_{d_i}) $ where $D$ is the 2D differential operator while the potential function $\phi(z) = \frac{|z|^p}{1 + |\frac{z}{c}|^{p-q}}$ with recommended choice $p = 2, q =1.5, c=10$; the second regularizer $g(x)=\iota_C(x)$ to be the indicator function of a box constraint restricting the pixel values to sit in the range of $[0,1]$. The image size is $256 \times 256$ and we consider sparse-view CT with $A\in\mathbb{R}^{46080\times65536}$ with 180 views. For data-fit we choose least-squares loss which is smooth, and also KL divergence which is nonsmooth. We consider the Condat-Vu's three-operator splitting scheme \citep{condat2013primal,vu2013splitting, chambolle2016ergodic} as the baseline -- this algorithm can be view as the deterministic counterpart of our algorithm. In both experiments we found that our algorithm converges significantly faster than the Condat-Vu's method in this type of three-composite problems with at least $O(1/K)$ rates, complimenting our theory. We then test our RED extensions of TOS-SPDHG in low-dose X-ray CT where the measurement noise is significant. Here we choose pretrained DnCNN as the denoiser for TOS-SPDHG-RED and TOS-SPDHG-eRED. We observe from Figure \ref{F2} and \ref{F3} that both of our RED extensions perform well in this task, while the TOS-SPDHG-eRED consistently perform the best in our experiments.

\section{Conclusion}

In this work we propose a stochastic primal-dual three operator splitting algorithm TOS-SPDHG admitting preconditioning and arbitrary sampling, for efficiently solving convex three-composite optimization problems, particularly in imaging applications. We build on the theoretical framework of \citep{chambolle2018stochastic} and extend the convergence proof to the case of three composites, showing the ergodic convergence rate of $O(1/K)$. Moreover, if any part of the function, either $f^*$, $g$ or $h$ is strongly-convex, we can trivially extend the current theory to prove faster convergence rates, which we omit here for now. Our algorithm is a simple yet crucial extension of SPDHG, with a much wider range of applications. Meanwhile, we show that the TOS-SPDHG can easily leverage deep networks for improved reconstruction with the RED framework.

\section*{Acknowledgments}

MJE acknowledges support from the EPSRC grants (EP/S026045/1, EP/T026693/1,
EP/V026259/1, EP/Y037286/1) and the European Union Horizon 2020 research
and innovation programme under the Marie Skodowska-Curie grant agreement REMODEL. CBS acknowledges support from the Philip Leverhulme
Prize, the Royal Society Wolfson Fellowship, the EPSRC advanced career fellowship EP/V029428/1, EPSRC
Grants EP/S026045/1 and EP/T003553/1, EP/N014588/1, EP/T017961/1, the Wellcome Innovator Awards
215733/Z/19/Z and 221633/Z/20/Z, the European Union Horizon 2020 research and innovation programme
under the Marie Skłodowska-Curie Grant agreement No. 777826 NoMADS, the Cantab Capital Institute for
the Mathematics of Information and the Alan Turing Institute.
\bibliography{main.bib}

\end{document}